**Generalized hypergeometric solutions of the Heun equation**

**A.M. Ishkhanyan**[1,2]

[1]Russian-Armenian University, Yerevan, 0051 Armenia
[2]Institute for Physical Research, NAS of Armenia, Ashtarak, 0203 Armenia

We present infinitely many solutions of the general Heun equation in terms of generalized hypergeometric functions. Each solution assumes that two restrictions are imposed on the involved parameters: a characteristic exponent of one of the singularities should be a non-zero integer, and the accessory parameter must satisfy a polynomial equation.

**PACS numbers:** 02.30.Gp Special functions, 02.30.Hq Ordinary differential equations, 02.30.Mv Approximations and expansions

**MSC numbers**: 33E30 Other functions coming from differential, difference and integral equations, 34B30 Special equations (Mathieu, Hill, Bessel, etc.), 30Bxx Series expansions

**Keywords:** general Heun equation; generalized hypergeometric function; recurrence relation

## 1. Introduction

The general Heun equation [1]

$$\frac{d^2u}{dz^2}+\left(\frac{\gamma}{z}+\frac{\delta}{z-1}+\frac{\varepsilon}{z-a}\right)\frac{du}{dz}+\frac{\alpha\beta z-q}{z(z-1)(z-a)}u=0\,, \tag{1}$$

where $$1+\alpha+\beta=\gamma+\delta+\varepsilon\,,$$

is widely used in contemporary fundamental and applied research (see, e.g., [2-6] and references therein). Nevertheless, this equation is much less studied than its immediate predecessor (the Gauss hypergeometric equation), and its solutions in terms of simpler functions, including functions of the hypergeometric class, are very rare.

Here, we present infinitely many solutions in terms of a single generalized hypergeometric function ${}_rF_s$ [7]. Our result is that such solutions exist if a characteristic exponent of one of the singularities is a nonzero integer and the accessory parameter $q$ satisfies a polynomial equation.

To be specific, we consider the singularity $z=a$. The characteristic exponents are $\mu_{1,2}=0,\,1-\varepsilon$. Let the exponent $\mu_2=1-\varepsilon$ be a non-zero integer. Our basic assertion is that for any *negative* integer $\varepsilon=-N$, $N=1,2,3,...$ (we discuss the case of a positive integer $\varepsilon$ later), the Heun equation admits the generalized hypergeometric solution

$$u={}_{N+2}F_{1+N}(1+e_1,...,1+e_N,\alpha,\beta;e_1,...,e_N,\gamma;z)\,. \tag{2}$$

This solution applies for certain particular choices of the accessory parameter $q$ defined by a polynomial equation of degree $N+1$. We note that for $\varepsilon=0$ the Heun equation admits a solution in terms of the ordinary hypergeometric function:

$N=0$: $\varepsilon=0$:

$$u={}_2F_1(\alpha,\beta;\gamma;z) \tag{3}$$

holds for

$$q-a\alpha\beta=0\,. \tag{4}$$

We present the solutions for $N=1$ and $N=2$.

$N=1$: $\varepsilon=-1$:

$$u={}_3F_2(\alpha,\beta,1+e_1;\gamma,e_1;z)\,, \tag{5}$$

$$(q-a\alpha\beta+a(1-\delta))(q-a\alpha\beta+(a-1)(1-\gamma))-a(1-a)(1+\alpha-\gamma)(1+\beta-\gamma)=0\,, \tag{6}$$

where the parameter $e_1$ is given by $e_1=a\alpha\beta/(q-a\alpha\beta)$. This quantity parameterizes the root of equation (6) as

$$q=a\alpha\beta\frac{1+e_1}{e_1}\,,\quad a=\frac{e_1\left(1+e_1-\gamma\right)}{\left(e_1-\alpha\right)\left(e_1-\beta\right)}\,, \tag{7}$$

This solution was first noticed by Letessier [8] and further studied by Letessier, Valent and Wimp [9] and by Maier [10].

$N=2$: $\varepsilon=-2$:

$$u={}_4F_3(\alpha,\beta,1+e_1,1+e_2;\gamma,e_1,e_2;z)\,, \tag{8}$$

$$\begin{gathered}\Big(\left(q-a\alpha\beta\right)^2+\left(q-a\alpha\beta\right)\left(\gamma-2+\left(1-\alpha-\beta\right)a\right)+2a\left(a-1\right)\alpha\beta\Big)\times\\ \left(q-a\alpha\beta-2\left(1+\alpha+\beta\right)a-2+2\gamma\right)+\left(q-a\alpha\beta\right)2a\left(a-1\right)\left(1+\alpha\right)\left(1+\beta\right)=0,\end{gathered} \tag{9}$$

where the parameters $e_{1,2}$ are the roots of the quadratic equation

$$e^2+\frac{3-2\gamma-q+a\left(1+2\alpha+2\beta+\alpha\beta\right)}{1-a}e-\frac{a\alpha\beta\left(2-2\gamma-q+a\left(2+2\alpha+2\beta+\alpha\beta\right)\right)}{\left(1-a\right)\left(q-a\alpha\beta\right)}=0\,. \tag{10}$$

This solution was obtained by Takemura [11]. It can be shown that these roots lead to the following relations resembling equations (7):

$$q=a\alpha\beta\frac{\left(1+e_1\right)\left(1+e_2\right)}{e_1e_2}\,,\quad a^2=\frac{e_1\left(1+e_1-\gamma\right)}{\left(e_1-\alpha\right)\left(e_1-\beta\right)}\frac{e_2\left(1+e_2-\gamma\right)}{\left(e_2-\alpha\right)\left(e_2-\beta\right)}\,, \tag{11}$$

In the general case $\varepsilon=-N$ the accessory parameter $q$ and the parameters $e_{1,2,\ldots,N}$ involved in solution (2) are determined from a system of $N+1$ algebraic equations. These equations are constructed by equating the coefficients of the following polynomial $\Pi(n)$ in an auxiliary variable $n$ to zero:

$$\Pi = a(\alpha-1+n)(\beta-1+n)\prod_{k=1}^{N}(e_k+n)+Q\prod_{k=1}^{N}(e_k-1+n)+(n-1)(\gamma-2+n)\prod_{k=1}^{N}(e_k-2+n), \quad (12)$$

where $$Q=-q+(n-1)(\delta+a\varepsilon)-(1+a)(n-1)(n-1+\alpha+\beta) \; .$$

A useful observation is that the polynomial $\Pi(n)$ has the degree $N$, not $N+2$ as might seen at first glance. This is because the two possible highest-degree terms proportional to $n^{N+1}$ and $n^{N+2}$ vanish. We therefore have $N+1$ equations, of which $N$ equations suffice for determining the parameters $e_{1,2,\ldots,N}$ and the remaining equation after elimination of $e_{1,2,\ldots,N}$ imposes a restriction on the parameters of the Heun equation. It can be verified that this restriction is a polynomial equation of degree $N+1$ for the accessory parameter $q$. For $N=0,1,2$, we obtain the corresponding equations (4), (6), and (9).

It can further be shown that for $q$ satisfying this equation, the parameters $e_{1,2,\ldots,N}$ can be determined as roots of a polynomial equation of degree $N$. For $N=1$ this equation is

$$e_1-\frac{a\alpha\beta}{q-a\alpha\beta}=0\,, \quad (13)$$

and the equation for $N=2$ is given by (10) (we discuss the equation for $N=3$ in the next section). We finally note that the system of algebraic equations for $e_{1,2,\ldots,N}$ leads to a generalization of equations (7) and (11):

$$q=a\alpha\beta\prod_{k=1}^{N}\frac{1+e_k}{e_k}\,, \qquad a^N=\prod_{k=1}^{N}\frac{e_k(1+e_k-\gamma)}{(e_k-\alpha)(e_k-\beta)}\,, \quad (14)$$

We discuss the derivation of the presented results in the next section.

## 2. Derivations

We consider the Frobenius series solution of general Heun equation (1) for the vicinity of the singularity $z=0$:

$$u=z^{\mu}\sum_{n=0}^{\infty}c_n z^n\,, \qquad \mu=0,\,1-\gamma\,. \quad (15)$$

The coefficients of this expansion satisfy the three-term recurrence relation

$$R_n c_n+Q_{n-1}c_{n-1}+P_{n-2}c_{n-2}=0\,. \quad (16)$$

For the exponent $\mu=0$ the coefficients of this relation are

$$R_n = a(\gamma - 1 + n)n ,$$

$$Q_n = -q + n(\delta + a\varepsilon) - (1 + a)n(n + \alpha + \beta) , \quad (17)$$

$$P_n = (\alpha + n)(\beta + n) .$$

The idea is to seek the cases where the Frobenius expansion (15) reduces to a generalized hypergeometric series. To study this possibility, we note that the generalized hypergeometric function ${}_rF_s$ is defined in terms of the series [7]

$${}_rF_s(a_1,...,a_r;b_1,...,b_s;z) = \sum_{n=0}^{\infty} c_n z^n , \quad (18)$$

where the coefficients obey the two-term recurrence relation

$$\frac{c_n}{c_{n-1}} = \frac{1}{n} \frac{\prod_{k=1}^{r}(a_k - 1 + n)}{\prod_{k=1}^{s}(b_k - 1 + n)} . \quad (19)$$

With function (2) in mind, we set $r = N + 2$, $s = N + 1$, and

$$a_1,....,a_N, a_{N+1}, a_{N+2} = 1 + e_1,...,1 + e_N, \alpha, \beta ,$$

$$b_1,....,b_N, b_{N+1} = e_1,...,e_N, \gamma . \quad (20)$$

We then rewrite recurrence relation (19) as

$$\frac{c_n}{c_{n-1}} = \frac{(\alpha - 1 + n)(\beta - 1 + n)}{(\gamma - 1 + n)n} \prod_{k=1}^{N} \frac{e_k + n}{e_k - 1 + n} . \quad (21)$$

Substituting this in equation (16), we obtain

$$R_n \frac{(\alpha - 1 + n)(\beta - 1 + n)}{(\gamma - 1 + n)n} \prod_{k=1}^{N} \frac{e_k + n}{e_k - 1 + n} + Q_{n-1} + P_{n-2} \frac{(\gamma - 2 + n)(n - 1)}{(\alpha - 2 + n)(\beta - 2 + n)} \prod_{k=1}^{N} \frac{e_k - 2 + n}{e_k - 1 + n} = 0 . \quad (22)$$

Substituting $R_n$ and $P_{n-2}$ here and canceling the common denominator, we rewrite this equation as

$$\Pi(n) = a(\alpha - 1 + n)(\beta - 1 + n) \prod_{k=1}^{N} (e_k + n) + Q_{n-1} \prod_{k=1}^{N} (e_k - 1 + n)$$
$$+ (n-1)(\gamma - 2 + n) \prod_{k=1}^{N} (e_k - 2 + n) = 0 . \quad (23)$$

This is a polynomial equation in $n$. Notably, the polynomial $\Pi(n)$ is of degree $N + 1$, not $N + 2$ because the possible highest-degree term $\sim n^{N+2}$ vanishes identically. Hence, equation (23) has the form

$$\Pi(n) = \sum_{m=0}^{N+1} A_m(a,q;\alpha,\beta,\gamma,\delta,\varepsilon;e_1,...,e_N) n^m = 0 . \quad (24)$$

Equating the coefficients $A_m$ to zero ensures the recurrence relation (16) is satisfied for all $n$.

We hence have $N+2$ equations $A_m = 0$, $m = 0,1,..,N+1$, of which $N$ equations serve to determine the parameters $e_{1,2,...,N}$ and the remaining two impose restrictions on the parameters of the Heun equation. One of these restrictions is easily derived by calculating the coefficient $A_{N+1}$ of the term proportional to $n^{N+1}$. Using the Fuchsian condition $1+\alpha+\beta = \gamma+\delta+\varepsilon$, we easily show that this coefficient is $(a-1)(\varepsilon+N)$. Hence, for $a \neq 1$,

$$\varepsilon = -N . \tag{25}$$

The second restriction is derived by eliminating $e_{1,2,...,N}$. For $N = 0,1,2$, these restrictions are given by the corresponding equations (4), (6), and (9). For higher $N$, the equations are cumbersome, and we omit them. But we note that this restriction can be alternatively derived via terminating the series solution of the Heun equation in terms of the Gauss hypergeometric functions [12]. This assertion follows because the generalized hypergeometric function (2) with $N$ numerator parameters exceeding the denominator parameters by unity has a representation as a linear combination with constant coefficients of a finite number of Gauss hypergeometric functions. This linear combination can be conveniently derived by truncating the expansions of the solutions of the Heun equation in terms of hypergeometric functions [12-15]. The termination condition for $\varepsilon = -N$ is a polynomial equation of degree $N+1$ for the accessory parameter $q$ [12].

Examining the structure of equations $A_m = 0$, we note that if $q$ is separately fixed by the abovementioned $(N+1)$-th degree equation, then we can write the equations $A_m = 0$ as $N$ linear equations for $N$ auxiliary variables $x_k$ given as $x_k = \sum_{1 \le i_1 < i_2 < \cdots < i_k \le N} e_{i_1} e_{i_2} \cdots e_{i_k}$, $k = 1,2,...,N$. It then follows that $e_{1,2,...,N}$ can be determined as roots of a polynomial equation of degree $N$. For $N = 1$ and $N = 2$, we have the respective the equations (13) and (10). For $N = 3$, the relations defining the polynomial equation for $e_{1,2,3}$ by Viet's theorem are

$$\begin{aligned} & e_1+e_2+e_3 = -\frac{6-3\gamma-q+a\left(3+3\alpha+3\beta+\alpha\beta\right)}{1-a}, \\ & e_1e_2+e_1e_3+e_2e_3 = \frac{22-12\gamma-3q+2a\left(4+3\alpha+3\beta\right)}{2\left(1-a\right)} - \\ & \frac{-2\gamma-q+a\left(2+\alpha\right)\left(2+\beta\right)}{2\left(1-a\right)}\left(e_1+e_2+e_3\right), \end{aligned} \tag{26}$$

$$e_1e_2e_3 = \frac{a\alpha\beta}{q-a\alpha\beta}\left(1+\left(e_1+e_2+e_3\right)+\left(e_1e_2+e_1e_3+e_2e_3\right)\right).$$

The parameters $e_{1,2,3}$ given by these equations provide a solution

$$u = {}_5F_4(\alpha,\beta,1+e_1,1+e_2,1+e_3;\gamma,e_1,e_2,e_3;z) \tag{27}$$

of the Heun equation if $\varepsilon=-3$ and $q$ satisfies the quartic equation

$$\begin{aligned}
&3a\left(\alpha-\gamma+1\right)\left(\beta-\gamma+1\right)\Big[3a(a-1)\left(\alpha-\gamma+3\right)\left(\beta-\gamma+3\right)+\\
&\left(3a\left(1-\delta\right)+X\right)\left(a\left(7-\gamma-2\delta\right)+\gamma-3+X\right)\Big]+\left(3\left(a-1\right)\left(1-\gamma\right)+X\right)\\
&\Big[4a\left(\alpha-\gamma+2\right)\left(\beta-\gamma+2\right)\left(3a\left(1-\delta\right)+X\right)+\left(a\left(7-2\gamma-\delta\right)+2\gamma-4+X\right)\\
&\left(3a\left(\alpha-\gamma+3\right)\left(\beta-\gamma+3\right)+\frac{\left(3a\left(1-\delta\right)+X\right)\left(a\left(7-\gamma-2\delta\right)+\gamma-3+X\right)}{a-1}\right)\Big],
\end{aligned} \tag{28}$$

where $X=q-a\alpha\beta$.

The derivation of the system of equations $A_m=0$, $m=0,1,..,N$, for the parameters $e_{1,2,...,N}$ completes the calculations. We additionally note that the solution of this system is unique up to an obvious permutation of $e_{1,2,...,N}$. Regarding the two equations (14), which can be useful for applications, we note that the first one for $q$ is a direct consequence of the equation $A_1=0$ for the coefficient of the term proportional to $n$. The second equation, which is a representation of the parameter $a$ through $e_{1,...,N}$, is easily derived by resolving the equations $\Pi(n=1-e_k)=0$, $k=1,2,...,N$, for $a$. For example, for $N=2$, the equations $\Pi(1-e_1)=0$ and $\Pi(1-e_2)=0$ give

$$a=-\frac{e_1\left(1+e_1-\gamma\right)\left(1+e_1-e_2\right)}{\left(e_1-\alpha\right)\left(e_1-\beta\right)\left(1-e_1+e_2\right)} \tag{29}$$

and

$$a=-\frac{e_2\left(1+e_2-\gamma\right)\left(1-e_1+e_2\right)}{\left(e_2-\alpha\right)\left(e_2-\beta\right)\left(1+e_1-e_2\right)}, \tag{30}$$

which immediately produce the second equation (11). We finally note that the second equation (14), although useful, is just a by-product. It is not needed for constructing the solution of the Heun equation. The $N+1$ equations $A_m=0$ suffice.

Let $\varepsilon$ now be a *positive* integer: $\varepsilon=N$, $N=1,2,3,...$. This case is easily treated by applying the elementary power-law change $u=(z-a)^{1-\varepsilon}w$, which transforms the Heun equation into another Heun equation with the altered parameter $\varepsilon_1=2-\varepsilon$. Indeed, for $\varepsilon\geq 2$, we have a Heun equation with a zero or negative integer $\varepsilon_1$. As a result, we obtain the solution

$$u=(z-a)^{1-\varepsilon}\ {}_{N+2}F_{1+N}(\tilde{e}_1+1,...,\tilde{e}_N+1,\alpha+1-\varepsilon,\beta+1-\varepsilon;\tilde{e}_1,...,\tilde{e}_N,\gamma;z)\,. \tag{31}$$

Thus, the only exception is the case $\varepsilon=1$ for which both characteristic exponents $\mu_{1,2}=0,1-\varepsilon$ are zero. We do not know a ${}_rF_s$ solution in this exceptional case.

Concluding this section, we note that for any integer $\varepsilon=N\in\mathbb{Z}$, $N\neq 1$, we can construct another set of similar solutions in terms of generalized hypergeometric functions by examining the Frobenius series solution of the Heun equation in the vicinity of the singularity $z=1$. The resulting solution for a negative integer $\varepsilon=-1,-2,...$ has the form

$$u={}_{N+2}F_{1+N}(e_1+1,...,e_N+1,\alpha,\beta;e_1,...,e_N,\delta;1-z)\,, \tag{32}$$

(we note that instead of $\gamma$, we here have $\delta$ as a denominator parameter), which in general is independent of solution (2). This solution hence represents the second independent fundamental solution of the Heun equation. Finally, for a positive integer $\varepsilon\neq 1$, we obtain a second independent fundamental solution as

$$u=(z-a)^{1-\varepsilon}\ {}_{N+2}F_{1+N}(e_1+1,...,e_N+1,\alpha+1-\varepsilon,\beta+1-\varepsilon;e_1,...,e_N,\delta;1-z)\,. \tag{33}$$

We note that for any set of the parameters, the generalized hypergeometric series involved in these solutions converge inside the unit circle centered at the corresponding singularity.

## 3. A physical example

The general Heun equation is currently encountered in many areaes of physics research ranging from classical physics and quantum mechanics to astronomy and cosmology. We present an example from quantum physics discussed in [16].

We consider the singular version of the third hypergeometric exactly solvable quantum-mechanical potential [17]

$$V=V_0+\frac{V_1}{\sqrt{1-e^{-x/\sigma}}}\,. \tag{34}$$

The one-dimensional Schrödinger equation

$$\frac{d^2\psi}{dx^2}+\frac{2m}{\hbar^2}\left(E-V(x)\right)\psi=0$$

for this potential reduces to Heun equation (1) with $\varepsilon=-1$ [18]. It can be verified that the parameters of this equation satisfy the equation (6) for $\varepsilon=-1$. As a result, we obtain the general solution of the Schrödinger equation in terms of the Clausen generalized hypergeometric functions ${}_3F_2$, which is written explicitly as

$$\psi = z^{\alpha_1}(z-1)^{\alpha_2}\left(c_1 u_1 + c_2 u_2\right) \tag{35}$$

with the two independent fundamental solutions

$$u_1 = {}_3F_2\left(\alpha,\beta,1+\frac{\alpha\beta}{q};\frac{\alpha\beta}{q},\gamma;z\right), \quad u_2 = {}_3F_2\left(\alpha,\beta,1-\frac{\alpha\beta}{q};-\frac{\alpha\beta}{q},\delta;1-z\right), \tag{36}$$

where $\gamma = 1+2\alpha_1$, $\delta = 1+2\alpha_2$, $q = \alpha_2 - \alpha_1$,

$$\alpha,\beta = \alpha_1 + \alpha_2 \pm \sqrt{\frac{8m\sigma^2}{\hbar^2}\left(-E+V_0\right)}, \tag{37}$$

$$\alpha_1 = \sqrt{\frac{2m\sigma^2}{\hbar^2}\left(-E+V_0-V_1\right)}, \quad \alpha_2 = \sqrt{\frac{2m\sigma^2}{\hbar^2}\left(-E+V_0+V_1\right)} \tag{38}$$

and $z = (1+\sqrt{1-e^{-x/\sigma}})/2$.

For a positive $\sigma$, potential (34) is defined on the positive half-axis $x > 0$. For $V_1 < 0$, the potential presents a bottomless well, which vanishes at infinity if $V_0 = -V_1$. Because it is a short-range potential, it supports only a finite number of bound states. These states are derived by requiring the wave function vanish both in the origin and infinity. The second of these requirements results in $c_1 = 0$, while the first provides the exact equation for the energy spectrum:

$${}_3F_2\left(\alpha,\beta,1-\frac{\alpha\beta}{q};-\frac{\alpha\beta}{q},\delta;\frac{1}{2}\right) = 0. \tag{39}$$

The exact number of bound states is equal to the number of zeros (not counting $x = 0$) of the zero-energy solution, which vanishes at the origin [19,20]. Because for $E = 0$ the exponent $\alpha_2$ vanishes and the lower parameter $\delta$ of the second independent solution (36) becomes unity, this solution should be constructed in a different way. The result is that the general solution of the Schrödinger equation for $E = 0$ is given by

$$\begin{aligned}\psi_{E=0} = c_1 z^{\alpha_1}\,{}_3F_2\left(\sqrt{2}\alpha_1+\alpha_1,-\sqrt{2}\alpha_1+\alpha_1,1+\alpha_1;\alpha_1,1+2\alpha_1;z\right)+\\ c_2 z^{-\alpha_1}\,{}_3F_2\left(\sqrt{2}\alpha_1-\alpha_1,-\sqrt{2}\alpha_1-\alpha_1,1+\alpha_1;\alpha_1,1;1-z\right).\end{aligned} \tag{40}$$

The condition $\psi_{E=0}(0) = 0$ then gives a linear relation between $c_1$ and $c_2$, which completes the construction of the zero-energy solution that vanishes in the origin.

## 4. Discussion

We have shown that the general Heun equation admits infinitely many particular solutions in terms of generalized hypergeometric functions which are simpler and significantly more studied functions. The existence of such solutions can be studied using the result obtained by Letessier, Valent and Wimp [9], who showed that the function ${}_{N+r}F_{s+N}$ with the structure of solution (2) satisfies a linear differential equation with polynomial coefficients of the order $\max(r,s+1)$. Their existence was conjectured by Takemura [11], who also proved this conjecture for $N \le 5$.

Compared with the known Heun-to-hypergeometric reductions constructed using a one-term ansatz involving a single Gauss hypergeometric function [21-24], the presented solutions have the advantage that they present cases with fewer restrictions on the parameters involved in the Heun equation. more precisely, each presented solution is obtained by imposing just two restrictions, while the mentioned reductions to the Gauss hypergeometric functions, except in trivial cases, assume three or more restrictions.

For the presented solutions, one restriction is imposed on a characteristic exponent of a singularity of the Heun equation: the exponent should be a non-zero integer. The other restriction is for the accessory parameter of the equation: this parameter must satisfy a certain polynomial equation. A straightforward examination then reveals that the corresponding singularity is *apparent*, i.e., the logarithmic term is absent in the general solution, and the solution is therefore analytic at this point [3]. The Heun equation with such a singularity can be regarded as a specific extension of the Gauss hypergeometric equation called the *deformed* hypergeometric equation [3,25,26]. It is known that in the case of the general and confluent Heun equations apparent singularities can be generated (or removed) by differentiating (or integrating) the equations [27]. But this is not the case for the hypergeometric equation [27,28]. Hence, the extension of the latter equation to a deformed version is a rather essential step. In particular, this extension allows generating the six Painlevé equations from the Heun equations via an anti-quantization procedure [29]. It is then remarkable that the solution of the deformed hypergeometric equation is written in terms of generalized hypergeometric functions. Moreover, it has been shown that the analytic continuation of a corresponding contour integral allows calculating the monodromy group for this equation in explicit terms [25,26].

We stress that each of the presented solutions can be written alternatively as a linear combination with constant coefficients of a finite number of Gauss hypergeometric functions [8-12,30-32]. These combinations are conveniently derived by truncating the series expansions of the solutions of the Heun equation in terms of the Gauss functions.

Regarding the applications of the presented solutions, we note that the two discussed restrictions imposed on the parameters of the Heun equation are satisfied in many physical situations. An example is the third exactly solvable hypergeometric quantum-mechanical potential [16], whose singular version we discussed here. Other examples include certain free-boundary problems [33] (e.g., coming from solidifcation [34] or fluid filtration [35] physics), gas dynamics [36], quantum mechanics [18,37], quantum two-state models [38], lattice systems [39], black hole physics [40], etc.

**Acknowledgments**

The work was supported by the Armenian State Committee of Science (SCS Grants Nos. 18RF-139 and 18T-1C276) and the Russian-Armenian (Slavonic) University at the expense of the Ministry of Education and Science of the Russian Federation.

**Conflict of Interest:** The authors declare that they have no conflicts of interest.